\documentclass [11pt,reqno]{amsart}
\usepackage {amsmath,amssymb,epsfig,enumerate,verbatim,geometry}

\newcommand{\C}{\mathbf{C}}

\renewcommand{\P}{\mathbf{P}}

\newcommand{\cB}{\mathcal{B}}

\newcommand{\cP}{\mathcal{P}}

\newcommand{\eg}{{\rm e.g.\ }} 
\newcommand{\ie}{{\rm i.e.\ }}

\DeclareMathOperator{\Exc}{Exc}

\numberwithin{equation}{section}       

\newtheorem{prop} {Proposition} [section]

\newtheorem{lem}[prop] {Lemma}

\newtheorem{prop-def}[prop]{Proposition-Definition}

\newtheorem*{mainthm}{Main Theorem} 
\theoremstyle{remark}

\newtheorem*{ackn}{Acknowledgment} 

\title{Endomorphisms of the plane preserving a pencil of curves}
\date{\today}

\author{Marius Dabija
  \and
  Mattias Jonsson}
\address{Dept of Mathematics\\
  University of Michigan\\
  Ann Arbor, MI 48109-1109\\
  USA\\
  and\\
  Dept of Mathematics\\
  KTH\\
  SE-100 44 Stockholm\\
  Sweden}

\email{mattiasj@umich.edu, mattiasj@kth.se} 
\thanks{Marius Dabija passed away on June 22, 2003.} 

\begin{document}

\begin{abstract}
  We classify endomorphisms of the plane that preserve a pencil of curves.
\end{abstract}

\maketitle


%
%
%
%
\section*{Introduction}
In dynamical systems, an interesting role is played by 
systems that are \emph{integrable} in the sense that they
preserve some additional structure. 
In this note we are concerned with
endomorphisms of the complex projective plane $\P^2$ that
preserve a pencil of 
curves.\footnote{See~\cite{Sibony} for a survey of dynamics
on projective spaces and~\cite{polyskew,fibered} 
for dynamics of endomorphisms appearing in the Main Theorem.}

More precisely, let $f:\P^2\to\P^2$ be an endomorphism
of algebraic degree $d>1$ and let $\cP$ be a pencil of
curves of degree $k\ge1$ on $\P^2$, defined by
$\pi:\P^2\dashrightarrow\P^1$.
For simplicity we shall assume that $\cP$ 
is \emph{irreducible}, that is, the generic curve in
$\cP$ is irreducible. (The reducible case is
discussed in Section~\ref{sec:extensions}.) We say that
$f$ leaves $\cP$ \emph{invariant} if $f$ maps any curve in $\cP$
(on)to another curve in $\cP$. This implies
that there exists an endomorphism $g:\P^1\to\P^1$ of degree $d$
such that $\pi\circ f=g\circ\pi$.

\begin{mainthm}
  There are only two types of invariant irreducible pencils $\cP$.
  \begin{itemize}
  \item[(i)]
    $\cP$ is an \emph{elementary} pencil, that is, the pencil of
    lines through a point. Pick homogeneous coordinates on $\P^2$ 
    and $\P^1$ such that 
    $\pi[x:y:z]=[x:y]$.
    Then $f$ takes the form $f[x:y:z]=[P(x,y):Q(x,y):R(x,y,z)]$.
  \item[(ii)]
    $\cP$ is a \emph{binomial} pencil. This means that we may
    pick homogeneous coordinates on $\P^2$ and $\P^1$ such that 
    $\pi[x:y:z]=[x^hy^{k-h}:z^k]$ where $0<h<k$ and $\gcd(h,k)=1$.
    In this case, after conjugation with an element of 
    $\mathrm{PGL}(3,\C)$ preserving $\cP$, $f$ takes the form
    $f[x:y:z]=[x^d:y^d:R]$ or 
    $f[x:y:z]=[y^d:x^d:R]$ for
    $R(x,y,z)=z^{d-kl}\prod_{i=1}^l(z^k+c_ix^hy^{k-h})$,
    where $0\le l\le d/k$ and $c_i\in\C^*$.
  \end{itemize}
\end{mainthm}

We shall only prove that an invariant (irreducible)
pencil $\cP$ must be elementary or binomial. The verification 
that $f$ is of the form stated in the theorem is straightforward
and left to the reader.

Our main theorem should be compared to the classification of curves
in $\P^2$ that are totally invariant by an endomorphism. Such
a curve must be a union of at most three 
lines~\cite{FS1,dabijathesis,CerveauLinsNeto,BCS},
a result that we exploit in the proof.

The strategy is to study the interplay between the
ramification loci of $\pi$, $f$ and $g$. In doing so,
the reducible and multiple curves in $\cP$ play an important role. 
We also use the minimal desingularization of $\cP$ and
the Riemann-Hurwitz formula to analyze the induced map between fibers.

%
%
%
%
\section{Counting ramification degrees}
Let $\cP$ be an irreducible pencil of plane curves of degree $k$
with associated projection $\pi:\P^2\dashrightarrow\P^1$, say
$\pi[x:y:z]=[A(x,y,z):B(x,y,z)]$ in homogeneous coordinates. 
By an \emph{element} of $\cP$ we shall mean the divisor of $sA+tB$
for some $[s:t]\in\P^1$. A \emph{curve} in $\cP$ will mean the 
(reduced) support of an element of $\cP$.

The \emph{base locus} $\cB(\cP)$ is the (nonempty, finite) 
indeterminacy set of $\pi$. The \emph{ramification locus} $R_\pi$ 
is defined as follows. For an effective divisor $D=\sum_jm_jD_j$ 
on $\P^2$, set $R_D=\sum_j(m_j-1)D_j$. Then $R_\pi=\sum R_D$, the 
sum being over the elements of $\cP$.
Set $e(\cP):=2k-\deg R_\pi$.
\begin{lem}\label{L1}
  We have $e(\cP)\ge 2$ with equality
  iff $\cP$ is an elementary pencil.
\end{lem}
\begin{proof}
  Consider a generic line $L\subset\P^2$. 
  Then $L\cap\cB(\cP)=\emptyset$, $L$ is
  transversal to $R_\pi$ and all tangencies (if any) between $L$
  and the curves in $\cP$ are quadratic.
  Applying the Riemann-Hurwitz formula to $\pi|_L:L\to\P^1$, 
  which has topological degree $k$, we obtain
  $2=2k-\deg R_\pi-t$, where $t\ge0$ is the number of 
  tangencies as above. 
  In particular, $t=0$ iff there are no tangencies and 
  this is easily seen to happen iff $\cP$ is elementary.
  Since $e(\cP)=2+t$, the proof is complete.
\end{proof}

Now assume that $f:\P^2\to\P^2$ is an endomorphism preserving
$\cP$. Let $g:\P^1\to\P^1$ be the induced endomorphism, that is,
$\pi\circ f= g\circ\pi$. 
Denote the ramification loci of $f$ and
$g$ by $R_f\subset\P^2$ and $R_g\subset\P^1$, respectively.
Write $R_f^\cP$ for the part of $R_f$ supported on curves in $\cP$.
A direct calculation based on $\pi\circ f= g\circ\pi$
yields
\begin{lem}\label{L3}
  We have $R_f^\cP+f^*R_\pi=R_\pi+\pi^*R_g$. 
  In particular, 
  $\deg R_f^\cP=e(d-1)$, with $e=e(\cP)$ as above.
\end{lem}
Since $\deg R_f=3(d-1)$ and $e\ge2$ we see that $e=2$ or $e=3$.
If $e=2$, then $\cP$ is an elementary pencil by Lemma~\ref{L1}, 
so let us assume from now on
that $e=3$. In particular, $R_f$ is supported on curves in $\cP$.
\begin{lem}\label{L4}
  Assume $\cP$ has a reducible curve.
  Then $\cP$ is a binomial pencil.
\end{lem}
\begin{proof}
  As $f$ maps curves in $\cP$ onto curves in $\cP$,
  the reducible curves in $\cP$ must be totally invariant for $f$.
  Since any totally invariant curve for $f$ is a union of at most
  three lines~\cite{FS1,dabijathesis,CerveauLinsNeto,BCS}, 
  there is exactly one
  reducible curve, the union of two or three lines.

  In the case of two lines, we assume $A(x,y,z)=x^hy^{k-h}$ with 
  $0<h<k$. 
  Let $m_jC_j$ be the multiple curves of $\cP$, \ie
  $m_j\ge2$ and $C_j$ is reduced and irreducible. Then
  $R_\pi=\sum_{j=1}^J(m_j-1)C_j+(h-1)\{x=0\}+(k-h-1)\{y=0\}$.
  Taking degrees and using $\deg R_\pi=2k-3$ we get
  $\sum_{j=1}^J(1-\frac{1}{m_j})=1-\frac1k$.
  Since $k>1$, the unique solution is given by 
  $J=1$, $m_1=k$.
  Thus $\cP$ admits a multiple line. Since $\cP$ is irreducible, 
  this line cannot pass through $[0:0:1]$, so we may assume
  $B(x,y,z)=z^k$. Then $\cP$ is a binomial pencil.

  When the reducible curve is a union of three lines we may assume 
  $A(x,y,z)=x^hy^lz^{k-h-l}$ with $h,l>0$, $h+l<k$. 
  Reasoning as above we get
  $\sum_{j=1}^J(1-\frac{1}{m_j})=1$,
  implying
  $J=2$ and $m_1=m_2=2$. But then $\cP$ is reducible,
  a contradiction.
\end{proof}
%
%
%
%
\section{Ruling out the remaining cases}
From now on we assume $e=e(\cP)=3$, that $\cP$ 
has no reducible curve, and that $f$ preserves $\cP$. 
In the end, we shall see that this is impossible. 

Let $p:X\to\P^2$ be a desingularization of the pencil $\cP$,
\ie $p$ is a composition of point blowups
and $\pi:\P^2\dashrightarrow\P^1$ lifts to a regular map
$\pi_X:X\to\P^1$ defining a ruling on $X$. 
Let $\Exc(p)\subset X$ be the
exceptional locus of $p$; we call its irreducible
components \emph{exceptional primes}.
There are finitely many exceptional primes
$E_i$, $i=1,\dots,I$ of $\Exc(p)$
that are \emph{horizontal}, \ie $\pi_X$ maps $E_i$
surjectively onto $\P^1$, say with degree $n_i$. Write
$E'$ for the union of the other (vertical) exceptional primes. 
Let $F_t=\pi_X^{-1}(t)$, $t\in\P^1$, be the fibers of 
$\pi_X$. For generic $t$, $F_t$ is smooth, 
intersects $E_i$ transversely at $n_i\ge1$ distinct generic points,
and is disjoint from $E'$. 

Recall that $R_f$ is supported on curves in $\cP$. 
Moreover, since $f$ is finite, 
the base locus $\cB(\cP)\subset\P^2$ is totally invariant. 
Thus, for generic $t$, $f$ induces a regular map of $F_t$ 
onto $F_{g(t)}$ and restricts to a covering map of 
$F_t\setminus\Exc(p)$ onto $F_{g(t)}\setminus\Exc(p)$.
By Riemann-Hurwitz this implies that $F_t$ is 
rational and that $\sum_{i=1}^I n_i=2$. 
In particular $\cP$ has one or two base points in $\P^2$.

As before, let $m_jC_j$ be the multiple curves of $\cP$.
Since $\cP$ has no reducible curves, the ramification locus
$R_\pi$ is supported on the $C_j$'s. 
Counting degrees as above, we get
$2-\frac3k=\sum_{j=1}^J(1-\frac1{m_j})$.
As $\cP$ is irreducible, the multiplicities $m_j$ are pairwise
relatively prime.
One easily checks that the only two possibilities are 
$k=2,J=1,m_1=2$ and 
$k=90$, $J=3$, $m_1=2$, $m_2=3$, $m_3=5$, respectively.

In the first case we have a pencil of conics with a
single double line. When there are two base points in $\P^2$,
one checks that $\cP$ has a reducible curve 
(and is binomial), a contradiction.
If there is a single base point in $\P^2$, we may assume
$A=y+x^2$, $B=y^2$. Then $I=1$, $n_1=1$ in
the notation above, contradicting $\sum n_i=2$.

In the second case, we have three multiple curves, 
of multiplicities 2, 3 and 5, respectively.
We may and will assume that the resolution 
$p:X\to\P^2$ of $\cP$ is \emph{minimal},  
a condition that may be paraphrased as follows. 
Let $\nu_i$ be the divisorial valuation on the function
field of $\P^2$ associated to the horizontal exceptional 
prime $E_i$, $1\le i\le I$.
Then $p:X\to\P^2$ is minimal with the property that
the centers of $\nu_i$ on $X$ are all one-dimensional.
By considering the dual graph of $p$ from this point of view
(see \eg \cite{spiv,valtree})
one finds that at least one of the horizontal 
exceptional primes, say $E_1$, intersects the
other exceptional primes in a finite set $Y$ consisting of
at most \emph{two} points. 
Thus at least one of the multiple fibers $F_t=\pi_X^{-1}(t)$, 
say $F_0$, must intersect $E_1$ but not $Y$.
As noted above, a generic fiber $F_t$ intersects $E_1$ transversely
at $n_1$ distinct points. Hence $(F_t\cdot E_1)=n_1$ and this equation
must then hold for \emph{all} $t$. In particular it holds for the multiple
fiber $F_0$. As $\sum n_i=2$, this is only possible if $I=1$, 
$n_1=2$ and $F_0=2G_0$ is
a double curve. Here $G_0$ is a priori reducible, but 
intersects $E_1$ transversely at a single point, not in $Y$.

Now pick a generic fiber $F_t$ and suppose that $g(t)=0$, that is,
the lift of $f$ to $X$ maps $F_t$ into $G_0$. 
Then the two 
intersection points of $F_t$ and $E_1$ must map to the 
unique intersection point of $G_0$ and $E_1$. 
Thus $G_0$ must be irreducible and we get a
covering map of degree $d>1$ from $F_t\setminus E_1$ onto 
$G_0\setminus E_1$. By Riemann-Hurwitz, this is impossible.

Hence, for $t$ generic we cannot have $g(t)=0$ and in fact not
$g^n(t)=0$ either for any $n\ge 1$. This implies that $t=0$ is
totally invariant for $g^2$, so the corresponding curve $C_0$
in $\P^2$ is totally invariant for $f^2$, a contradiction since 
$\deg C_0=45>3$.

This concludes the proof of the Main Theorem.
%
%
%
%
\section{Reducible pencils and further remarks}\label{sec:extensions}
Finally we discuss the case when $\cP$ is reducible (but without
base divisor). One must then be a little careful with the notion of
invariance. For example, let $\cP$ be defined by
$\pi[x:y:z]=[x^2:y^2]$ and set $f[x:y:z]=[x^2:y^2:z^2]$.
Then $f$ maps curves in $\cP$ into, but not onto, curves in $\cP$. 
On the other hand, the preimage of a fiber is
a union of fibers and we shall use this as the definition of
invariance.

It is a result of Bertini that $\cP$ factors through an
irreducible pencil $\cP'$, say defined by
$\pi':\P^2\dashrightarrow\P^1$.
This means that $\pi=\phi\circ\pi'$, where $\phi:\P^1\to\P^1$
is an endomorphism of degree $>1$. 
Clearly, if $\cP$ is invariant by $f$, then so is $\cP'$, which 
implies that $\cP'$ is elementary or binomial and that $f$ is
as described in the Main Theorem. 

To completely solve the classification problem in the reducible case
would therefore amount to the following problem. Given an
endomorphism $g'$ of $\P^1$ associated to one of the situations in
the Main Theorem, find all endomorphisms $\phi$ of $\P^1$ of 
degree~$>1$ such that $\phi\circ g'=g\circ\phi$ for another
endomorphism $g$. 
This appears to be an interesting and nontrivial problem,
but one of a different flavor than the present note, 
and we shall leave it for now.

Beyond the setup used here, there are numerous possible notions
of integrable holomorphic dynamics.
One could look at both more general (rational) mappings,  
more general structures (linear systems, webs, foliations,\dots)
and of course higher dimensions.
In~\cite{invwebs} we study algebraic webs on $\P^2$ invariant 
under endomorphisms.
Cantat and Favre~\cite{CF} have classified 
infinite groups of birational surface maps preserving a 
holomorphic foliation. 
Pereira~\cite{invfols} has outlined a generalization of the
results in this paper from pencils to foliations.
However, many interesting problems remain open.

\begin{ackn}
  The second author was supported by the NSF,
  the Swedish Research Council and the Gustafsson Foundation.
  He is grateful to Igor Dolgachev, Charles Favre and Jorge Pereira
  for valuable discussions. He also thanks 
  the referee for suggesting substantial 
  simplifications to several arguments in
  a previous version of the paper.
\end{ackn}
%
%
%
%

\end{document}